%% file: paynethesis.tex
\input hyperbasics

\input amstex
\documentstyle{phd}
\input btxmac
\overfullrule=0pt

\Author Garth Payne

\ThesisTitle Multivariate Hypergeometric Terms

\Date December 1997

\IncludeCopyright:yes
\Draft:no
\IncludeFrontMatter:yes

\loadmsam
\loadmsbm
\loadbold

\newsymbol\leqslant 1336 \let\le=\leqslant
\newsymbol\geqslant 133E 
\newsymbol\nil 203F
\newsymbol\smallfrown 1361
\newsymbol\nsubseteq 232A

\ifx \hyperref \undefined
  \def \hyperref#1#2{{#2}}
\fi
\ifx \hyperdef \undefined
  \def \hyperdef#1#2#3#4{{#4}}
\fi

\def\Statement[#1]{}

\def\Mark#1\StopMark{{#1}}
\def\OreSato{{\tt[Ore{\rm--}Sato]}}

\def\R{{\Cal R}}
\def\Q{{\Bbb Q}}
\def\Z{{\Bbb Z}}

\def\vec#1{{\boldkey#1}}

\def\gp#1#2#3{\mathop{\smash{\prod\limits_{#1}}\vphantom{\prod}}\nolimits_{#2}^{#3}}
\def\z{\vec z}
\def\e{\vec e}
\def\w{\vec w}
\predefine\hacek{\v}
\def\v{\vec v}
\def\u{\vec u}

\def\newop#1{\expandafter\def\csname #1\endcsname
        {\mathop{\text{\rm #1}}\nolimits}}
        
\newop{inv}
\newop{rgal}
\newop{fix}
\newop{relation}
\newop{sg}

\def\zd{zero divisor}
\def\hgt{hypergeometric term}
\def\c{{\vec c}}

\def\m1{^{-1}}

\Start

\Committee
\Member George E.~Andrews
Evan Pugh Professor of Mathematics
Head of the Department of Mathematics
Thesis Advisor
Chair of Committee
\Member Leonid Vaserstein
Professor of Mathematics
\Member W.~Dale Brownawell
Distinguished Professor of Mathematics
\Member Calyampudi R. Rao
Professor and Holder of the Eberly Family Chair in Statistics

\SignatoryPage

\input abs

\Contents

\Acknowledgements{ack}

\InputFile:one

\InputFile:two

\InputFile:three

\InputFile:four

\InputFile:five

\nocite{*}

\Bibliography{garth}

\InputFile:errata

\Stop

%% file: hyperbasics.tex
\expandafter\edef\csname hypers@fe\endcsname{\catcode
                                             `\noexpand @=\the\catcode`\@}%
\catcode`\@=11
%
%
\ifx\hyperd@ne\hyper@ndefined
 \global\let\hyperd@ne=\relax
\else
 \errhelp{hyperbasics.tex needs to be included only once outside
          of any {...} or \begingroup...\endgroup. You have tried to
          include it more than once. If the previous include was indeed
          outside any groupings, continue and all will be well.}%
 \errmessage{Input this file only once!}%
  
\fi
%
%
\def\hyperv@rsion{8}%
%
%
\newread\hyperf@le
\def\hyperf@lename{\jobname.hrf}%
\immediate\openin\hyperf@le\hyperf@lename\relax
\ifeof\hyperf@le\relax
 \immediate\closein\hyperf@le\relax
\else
 \immediate\closein\hyperf@le\relax
 \input \hyperf@lename
\fi
%
%
\newwrite\hyperf@le
\immediate\openout\hyperf@le\hyperf@lename
%
%
\newtoks\hypert@ks
%
%
\edef\hypert@mp{\catcode`\noexpand\#=\the\catcode`\#}%
\catcode`\#=12
\def\hyperh@sh{#}%
\hypert@mp
\let\hypert@mp=\relax
\let\hyper@nd=\relax
\def\hyperstr@pquote"#1"#2\hyper@nd{\ifx\hyper@ndefined#2\hyper@ndefined#1\else
                                    \ifx\hyper@ndefined#1\hyper@ndefined
                                    \hyperstr@pquote#2"\hyper@nd\else
                                    #1\hyperstr@pquote"#2"\hyper@nd\fi\fi}%
\def\hyperstr@pblank" #1 #2\hyper@nd"{\ifx\hyper@ndefined#2\hyper@ndefined#1\else
                                    \ifx\hyper@ndefined#1\hyper@ndefined
                                    \hyperstr@pblank"#2 \hyper@nd"\else
                                    #1\hyperstr@pblank" #2 \hyper@nd"\fi\fi}
\long\def\hyper@nchor#1#2{\edef\hyperm@cro{html:<A #1>}%
                          \special\expandafter{\hyperm@cro}%
                          {#2}}%
\def\hyper@atm@ning#1->#2\hyper@nd{#2}
\def\hyperlink#1{\edef\hypert@mp{#1}%
               \edef\hypert@mp{\expandafter\hyper@atm@ning\meaning\hypert@mp
                               \hyper@nd}%
               \edef\hypert@mp"{ \expandafter\hyperstr@pquote\expandafter"%
                               \hypert@mp"\hyper@nd}%
               \edef\hypert@mp{\expandafter\hyperstr@pblank\expandafter%
                               "\hypert@mp" \hyper@nd"}%
               \hyper@nchor{href=\expandafter"\hypert@mp"}}%
\def\hypertarget#1{\edef\hypert@mp{#1}%
               \edef\hypert@mp{\expandafter\hyper@atm@ning\meaning\hypert@mp
                               \hyper@nd}%
               \edef\hypert@mp"{ \expandafter\hyperstr@pquote\expandafter"%
                               \hypert@mp"\hyper@nd}%
               \edef\hypert@mp{\expandafter\hyperstr@pblank\expandafter%
                               "\hypert@mp" \hyper@nd"}%
               \hyper@nchor{name=\expandafter"\hypert@mp"}}%
\def\hyperref{\afterassignment\hyperr@f\let\hyperp@ram}
\def\hyperr@f{\ifx\hyperp@ram{\iffalse}\fi
               \expandafter\expandafter\expandafter\hyperr@@
               \expandafter{%
              \else
               \iffalse}\fi
               \ifx\hyperp@ram\hyper@ndefined
                 \message{Undefined reference}%
                 \def\hyperp@r@m{{}{undefined}{}}%
               \else
                 \edef\hyperp@r@m{\hyperp@ram}%
               \fi
               \expandafter\expandafter\expandafter\hyperr@@
               \expandafter\hyperp@r@m
              \fi}%
\def\hyperr@@#1#2#3{\ifx\hyper@ndefined#1\hyper@ndefined
                    \hypert@ks\expandafter{\hyperh@sh#2.#3}%
                    \else
                     \ifx\hyper@ndefined#2#3\hyper@ndefined
                      \hypert@ks{#1}%
                     \else
                      \def\hypert@mp{#1}%
                      \hypert@ks\expandafter\expandafter\expandafter
                      {\expandafter\hypert@mp\hyperh@sh#2.#3}%
                     \fi
                    \fi
                    \expandafter\hyperlink\expandafter{\the\hypert@ks}}%
\def\hyperdef#1#2#3{{\global\escapechar=`\\\relax
                     \edef\hypert@mp{\hyperstr@pquote"#2.#3"\hyper@nd}%
                     \expandafter\ifx\csname hyperd@\meaning\hypert@mp
                     \endcsname
                     \relax
                     \expandafter\gdef\csname hyperd@\meaning\hypert@mp
                     \endcsname{}%
                     \gdef#1{{}{\hyperstr@pquote"#2"\hyper@nd}%
                               {\hyperstr@pquote"#3"\hyper@nd}}%
                     \immediate\write\hyperf@le{\def\noexpand#1{#1}}%
                     \xdef\hypert@mp{\global\let\noexpand\hypert@mp=\relax
                                     \noexpand\hypertarget{\hypert@mp}}%
                     \global\hypert@ks={\hypert@mp}%
                     \else
                     \message\expandafter{'\hypert@mp' duplicate}%
                     \global\let\hypert@mp=\relax
                     \global\hypert@ks={\hyperdef{#1}{#2}{#3@}}%
                     \fi}\the\hypert@ks}%

\def\hyper@nique#1#2#3#4{\global\escapechar=`\\\relax
                     \edef\hypert@mp{\hyperstr@pquote"#2.#3"\hyper@nd}%
                     \expandafter\ifx\csname hyperd@\meaning\hypert@mp
                     \endcsname
                     \relax
                     \gdef#1{{}{\hyperstr@pquote"#2"\hyper@nd}%
                               {\hyperstr@pquote"#3"\hyper@nd}}%
                     \global\let\hypert@mp=\relax
                     #4%
                     \else
                     \global\let\hypert@mp=\relax
                     \hyper@nique{#1}{#2}{#3@}{#4}%
                     \fi
                     }%

\let\hyper@@@@=\relax
\def\hyper@@{\let\hyper@@@=\relax}%
\hyper@@
\def\hyper@{\relax\let\hyper@@@\noexpand\hyper@\noexpand}%
\def\hyperpr@ref{\hyper@@\hyperref}
\def\hyperpr@def{\hyper@@\hyperdef}

\let\href\hyperlink

%
%
\hypers@fe

%% file: abs.tex
In \hyperref\WZninetytwo{\bibref[WZ92]}, Wilf and Zeilberger developed a proof theory for
hypergeometric multisums centered around the notion of a multivariate
\hgt. A multivariate function $f(z_1,\ldots,z_k)=f(\z)$ from $\Z^k$
to a field $K$ is a {\it\hgt} if for each $i\in\{1,\ldots,k\}$ there
exist nonzero polynomials $A_i(\z)$ and $B_i(\z)$ in $K[\z]$ such
that
$$A_i(\z)f(\z)=B_i(\z)f(\z+\e_i)$$
for all $\z\in\Z^k$. Here $\e_1,\ldots,\e_k$ is the
standard basis for $\Z^k$. If $f$ is not a {\it\hyperref\defBtwo{\zd}}, then the {\it \hyperref\defBfour{term
ratios}} $R_1=A_1/B_1,\ldots,R_k=A_k/B_k$ are unique and satisfy the
relation
$$R_i(\z)R_j(\z+\e_i)=R_j(\z)R_i(\z+\e_j)\quad
\text{for each }i,j\in
\{1,\ldots,k\}$$
(cf.\ \hyperref\WZninety{\bibref[WZ90]}).

We introduce the concepts of divisibility lattice paths, \hyperref\defAthree{rational Galois  
spaces}, and \hyperref\defAfour{fixed factors} of rational functions to the study of the
relation for the \hyperref\defBfour{term ratios}. We prove that a solution $R_1,\ldots,R_k$
must be of the form \hyperref\oresato{\OreSato}
$$
R_i(\z)=\frac{C(\z+\e_i)}{C(\z)}\frac{D(
\z)}{D(\z+\e_i)}\prod_{\v\in V}\gp j0{\e_i \cdot \v}\frac{a_{
\v}(\v\cdot\z+j)}{b_{\v}(\v\cdot\z+j)}
\text{\ \ for $i\in\{1,\ldots,k\}$,}
$$
where $C$ and $D$ are polynomials, $V$ is a
finite subset of $\Z^k$, and, for each $\v\in V$, $a_{\v}$ and
$b_{\v}$ are {\sl univariate} polynomials all independent of
$i$.
The symbol $\gp jab{}$ denotes 
$\prod_{i=a}^{b-1}$ if $b>a$,
$1$ if $a=b$,
and
$1/\prod_{i=b}^{a-1}$ if $a>b$.

We use this factorization of $R_i$ to answer an obvious question about
multivariate \hgt s. Recall the Pochhammer symbol $(m)_r=m(m+1)\cdots(m+r-1)$. The
multivariate \hgt s that arise in practice have the form $$f(
\z)=\gamma_1^{z_1}\cdots\gamma_k^{z_k}\frac{C(\z)}{D(
\z)}\frac{\prod_{i=1}^p(m_i)_{\v_i\cdot
\z+r_i}}{\prod_{j=1}^q(n_j)_{\w_j\cdot\z+s_i}}$$ where
$\gamma_1,\ldots,\gamma_k\in K$, $C$ and $D$ are polynomials in $K[
\z]$, the $\v_i$ and $\w_j$ are in $\Z^k$, the $r_i$ and $s_j$
are in $\Z$, and the $m_i$ and $n_j$ are in $K$. The question is: do
all \hgt s have this form? We prove that if $K$ is algebraically
closed, and the \hgt\ is \hyperref\defBthirteena{honest}, then such an expression for $f$
exists \Mark\hyperref\Honestpw{{\it piecewise}}\StopMark. This is trivial in the case of one variable,
but {\sl not} in the case of several variables. We use this result to
settle the discrete part of a conjecture of Wilf and Zeilberger
\hyperref\WZninetytwo{\bibref[WZ92]} by showing that a holonomic \hgt\ is \hyperref\defpwp{piecewise proper},
which means roughly that we can take $D(\z)=1$ in the expression for 
$f$ above.

For any subsets $A$ and $B$ of an additive group $G$, define
$A+B=\{a+b:a \in A \text{ and } b \in B \}$ and $-A= \{-a:a
\in A \}$.  A subset $S$ of $G$ is said to be sum-free,
complete, and symmetric respectively if $S+S \subset S^\c$,
$S+S \supset S^\c$, and $S=-S$.  Cameron asked if for all
sufficiently large moduli $m$ there exists a sum-free complete
set in $\Z/m\Z$ that is not symmetric \hyperref\CamPortrait{\bibref[CamPortrait]}.
We answer Cameron's question by showing there
exists such a set for all moduli greater than or equal to
$890626$.  We also show that every sum-free complete set in
$\Z/m\Z$ that is not symmetric can be used to construct a
counterexample to a conjecture of Conway disproved by Marica
\hyperref\Marica{\bibref[Marica]}.  Conway conjectured that for any finite set
$S$ of integers, $|S+S| \le |S-S|$.

%% file: paynethesis.bbl
\hyperdef\chapteRBibliography{ChapterSection}{chapteRBibliography}{}
\begin{thebibliography}{10}

\bibitem{And76}
G.E. Andrews.
\newblock The theory of partitions.
\newblock In G.C. Rota, editor, {\em Encyclopedia of Mathematics and Its
  Applications}, volume~2. Addison-Wesley, 1976.%
\hyperdef\Andseventysix{Bibliography}{Andseventysix}{}

\bibitem{Ber1}
I.N. Bernstein.
\newblock Modules over a ring of differential operators, study of the
  fundamental solutions of equations with constant coefficients.
\newblock {\em Functional Analysis and Its Applications}, 5(2):89--101, 1971.%
\hyperdef\Berone{Bibliography}{Berone}{}

\bibitem{Bjo1}
J.~Bjork.
\newblock {\em Rings of Differential Operators}.
\newblock North-Holland, Amsterdam, 1979.%
\hyperdef\Bjoone{Bibliography}{Bjoone}{}

\bibitem{Calkin}
N.~Calkin.
\newblock On the number of sum-free sets.
\newblock {\em Bull. London Math. Soc.}, 22:141--144, 1990.%
\hyperdef\Calkin{Bibliography}{Calkin}{}

\bibitem{CamGraph}
P.~Cameron.
\newblock Cyclic automorphisms of a countable graph and random sum-free sets.
\newblock {\em Graphs and Combinatorics}, 1:129--135, 1985.%
\hyperdef\CamGraph{Bibliography}{CamGraph}{}

\bibitem{CamStructure}
P.~Cameron.
\newblock On the structure of a random sum-free set.
\newblock {\em Probability Theory and Related Fields}, 76:523--531, 1987.%
\hyperdef\CamStructure{Bibliography}{CamStructure}{}

\bibitem{CamPortrait}
P.~Cameron.
\newblock Portrait of a typical sum-free set.
\newblock In C.~Whitehead, editor, {\em Surveys in Combinatorics 1987}, pages
  13--42. Cambridge Univ. Press, Cambridge, 1987.%
\hyperdef\CamPortrait{Bibliography}{CamPortrait}{}

\bibitem{Gosp78}
R.W. Gosper.
\newblock Decision procedures for indefinite hypergeometric summation.
\newblock {\em Proc. Natl. Acad. Sci. U.S.A.}, 75:40--42, 1978.%
\hyperdef\Gospseventyeight{Bibliography}{Gospseventyeight}{}

\bibitem{Marica}
J.~Marica.
\newblock On a conjecture of {C}onway.
\newblock {\em Canad. Math. Bull.}, 12:233--234, 1969.%
\hyperdef\Marica{Bibliography}{Marica}{}

\bibitem{Pet92}
M.~Petkov\hacek{s}ek.
\newblock Hypergeometric solutions of linear recurrences with polynomial
  coefficients.
\newblock {\em J. Symbolic Computation}, 14:243--264, 1992.%
\hyperdef\Petninetytwo{Bibliography}{Petninetytwo}{}

\bibitem{Schur}
I.~Schur.
\newblock \"{U}ber die kongruenz $x^m+y^m\cong z^m \pmod p$.
\newblock {\em Jahresber. Deutsch. Math.-Verein.}, 25:114--117, 1916.%
\hyperdef\Schur{Bibliography}{Schur}{}

\bibitem{Stein}
S.~Stein.
\newblock The cardinalities of $a+a$ and $a-a$.
\newblock {\em Canad. Math. Bull.}, 16(3):343--345, 1973.%
\hyperdef\Stein{Bibliography}{Stein}{}

\bibitem{Wallis}
W.~Wallis, A.~Street, and J.~Wallis.
\newblock {\em Combinatorics: Room Squares, Sum-Free Squares, Hadamard
  Matrices}.
\newblock Springer-Verlag, Berlin, 1972.%
\hyperdef\Wallis{Bibliography}{Wallis}{}

\bibitem{WZ90}
H.S. Wilf and D.~Zeilberger.
\newblock Rational functions certify combinatorial identities.
\newblock {\em J. Amer. Math. Soc.}, 3:147--158, 1990.%
\hyperdef\WZninety{Bibliography}{WZninety}{}

\bibitem{WZ92}
H.S. Wilf and D.~Zeilberger.
\newblock An algorithmic proof theory for hypergeometric (ordinary and ``$q$'')
  multisum/integral identities.
\newblock {\em Invent. Math.}, 108:575--633, 1992.%
\hyperdef\WZninetytwo{Bibliography}{WZninetytwo}{}

\bibitem{Zeil90}
D.~Zeilberger.
\newblock A holonomic systems approach to special functions identities.
\newblock {\em J. of Computational and Applied Mathematics}, 32:321--368, 1990.%
\hyperdef\Zeilninety{Bibliography}{Zeilninety}{}

\end{thebibliography}
